\documentclass[12pt]{article}
\def\finproof{\hfill\hbox{\vrule width1.0ex height1.5ex}\vspace{2mm}}

\begin{document}

\begin{center}
{\Large 
Differential equations with discrete state-dependent delay:
uniqueness and well-posedness in the space of continuous
functions}

\bigskip

{\sc Alexander V. Rezounenko}

\smallskip

Department of Mechanics and Mathematics, Kharkov
University,

 4, Svobody Sqr., Kharkov, 61077, Ukraine

\smallskip

 E-mail: rezounenko@univer.kharkov.ua

\end{center}

\begin{quote}
{\bf Abstract.} Partial differential equations with discrete
(concentrated) state-dependent delays in the space of continuous
functions are investigated. In general, the corresponding initial
value problem is not well posed, so we find an additional
assumption on the state-dependent delay function to guarantee the
well posedness. For the constructed dynamical system we study the
long-time asymptotic behavior and prove the existence of  a
compact global attractor.
\end{quote}

\medskip



\bigskip

{\bf 1. Introduction} 

\marginpar{\tiny Feb 4, 2008}


\bigskip 

Delay differential equations is one of the oldest branches
of the theory of infinite dimensional dynamical systems -
theory which describes qualitative properties of systems,
changing in time.

We refer to the classical monographs on the theory of ordinary
(O.D.E.) delay equations
\cite{Hale,Hale_book,Walther_book,Azbelev,Mishkis}. The theory of
partial (P.D.E.) delay equations is essentially less developed
since such equations are infinite-dimensional in both time (as
delay equations) and space  (as P.D.E.s) variables, which makes
the analysis more difficult. We refer to some works which are
close to the present research \cite{travis_webb,
Chueshov-JSM-1992,Cras-1995,NA-1998,Rezounenko-2003} and to the
monograph \cite{Wu_book}.

Recently, a new class of delay equations - equations with a
state-dependent delay (SDD) attracts much attention of
researchers (see e.g. \cite{Walther2002,Walther_JDE-2003,
Nussbaum-Mallet-1992, Nussbaum-Mallet-1996, MalletParet,
Krisztin-2003, Walther_JDDE-2007} and also the review
\cite{Hartung-Krisztin-Walther-Wu-2006} for details and
references). Investigations of these equations essentially
differ from the ones of equations with constant or
time-dependent delays. The main difficulty is that
nonlinearities with SDD (in contrast to constant or
time-dependent delays) are not Lipschitz continuous on the
space of continuous functions - the main space of initial
data, where the classical theory of delay equations is
developed (see the references above). As a consequence, the
corresponding initial value problem (IVP), in general, is
{\bf not} well posed (in the sense of J.~Hadamard
\cite{Hadamard-1902,Hadamard-1932}).
 An explicit example of the nonuniqueness of solutions
to an ordinary equation with state-dependent delay (SDD) is given in
\cite{Winston-1970} (see also
\cite[p.465]{Hartung-Krisztin-Walther-Wu-2006}). As noticed in
\cite[p.465]{Hartung-Krisztin-Walther-Wu-2006} "typically, the IVP
is uniquely solved for initial and other data which satisfy suitable
Lipschitz conditions." The idea to investigate SDD O.D.E.s in the
space of Lipschitz continuous functions is very fruitful (see the
references above).

Unfortunately, in contrast to O.D.E.s with state-dependent delay,
one has {\bf no} possibility to exploit the space of Lipschitz
continuous functions for P.D.E.s with SD delay, because solutions
of P.D.E.s usually do not belong to this space.

The first attempt to treat P.D.E.s with state-dependent
(state-selective) delays has been made for a distributed delay
problem~\cite{Rezounenko-Wu-2006,Rezounenko_JMAA-2007}. The
existence (without uniqueness) of solutions for P.D.E.s with
discrete state-dependent delay was studied in \cite{Hernandez-2006}
(mild solutions) and \cite{Rezounenko_JMAA-2007} (weak solutions).

As noticed above, it is a common point of
view~\cite{Hartung-Krisztin-Walther-Wu-2006} that equations (O.D.E.s
and P.D.E.s) with discrete state-dependent delay are {\bf not} well
posed in the space of continuous functions (C). This leads to the
search of (particular) classes of equations which may be well-posed
in the space of continuous functions (C).

The main goal of the present paper is to propose an assumption on
the delay, which is sufficient for the well-posedness of the
corresponding initial value problem in the space C. To the best of
our knowledge, this is the first result for the well-posedness in C
of the discrete state-dependent delay (for P.D.E.s as well as for
O.D.E.s).

Discussing the meaning of the main assumption (H) (see below) for
applied problems, we hope that the assumption is the natural
mathematical expression of the fact that for many applied problems,
the models have a parameter (time $\eta_{ign}>0$) which is necessary
to take into considerations the time changes in the system. The
changes not always can be taken into considerations immediately. To
this end, the existence of $\eta_{ign}>0$ (no matter how small the
value of $\eta_{ign}>0$ is!) makes the corresponding initial value
problem well posed.

\bigskip {\bf 2. Formulation of the model with state-dependent discrete delay}

\medskip

Our  goal is to present an approach to study the following
 partial differential equation with {\it state-dependent
\underline{discrete} delay}

\begin{equation}\label{sdd8-1}
\frac{\partial }{\partial t}u(t,x)+Au(t,x)+du(t,x)
=  \int_\Omega b(u(t-\eta (u_t), y)) f(x-y) dy \equiv \big(
F(u_t)
\big)(x),\quad x\in \Omega ,
\end{equation}
 where $A$ is a densely-defined self-adjoint positive linear operator
 with domain $D(A)\subset L^2(\Omega )$ and with compact
  resolvent, so $A: D(A)\to L^2(\Omega )$ generates an analytic semigroup,
  $\Omega $ is a smooth bounded domain in $R^{n_0}$,
 $f: \Omega -\Omega \to R$   is a bounded function to
  be specified later,
$b:R\to R$ is a locally Lipschitz 
map, 
$d$ is a non-negative constant. The function
  $\eta (\cdot): C([-r,0];L^2(\Omega)) \to R_{+}$ represents the
state-dependent {\it discrete} delay.
We denote for short $C\equiv C([-r,0];L^2(\Omega)).$ The norms in
$L^2(\Omega)$ and $C$ are denoted by $||\cdot ||$ and $||\cdot ||_C$
respectively. As usually for delay equations, we denote by $u_t$ the
function of $\theta\in [-r,0]$ by the formula $u_t\equiv u_t(\theta)\equiv
u(t+\theta).$

We consider equation (\ref{sdd8-1}) with the following
initial condition
\begin{equation}\label{sdd8-ic}
u|_{[-r,0]}=\varphi \in C\equiv C([-r,0];L^2(\Omega)).
\end{equation}

The methods used in our work can be applied to another types of
nonlinear and delay P.D.E.s (as well as O.D.E.s). We choose a
particular form of nonlinear delay terms $F$ for simplicity and to
illustrate our approach on the diffusive Nicholson's blowflies
equation (see the end of the article for more details). Below you
will also find a remark (Rem.~11) on the local in space variable
problems.


\bigskip

{\bf 3. The existence of mild solutions 
}%

\medskip

In our study we use the standard

\smallskip

{\bf Definition~1}. 
 {\it A function $u\in C([-r,T]; L^2(\Omega))$ is called a {\tt mild solution}
 on $[-r,T]$ of the initial value problem (\ref{sdd8-1}), (\ref{sdd8-ic}) if it satisfies
 (\ref{sdd8-ic}) and
 \begin{equation}\label{sdd8-3-1}
u(t)=e^{-A t}\varphi(0) + \int^{t}_0 e^{- A (t-s)} \left\{ F(u_s) - d
\cdot u(s)\right\}\, ds, \quad t\in [0,T].
 \end{equation}
}%

\medskip

{\bf Proposition~1}. {\it Assume the function $b : R\to R$ is a locally
Lipschitz 
map, satisfying $|b(w)|\le C_1|w|+C_b$ with $C_i\ge 0,$ and
function
  $\eta (\cdot): C([-r,0];L^2(\Omega)) \to R_{+}$ is continuous,
  $f: \Omega -\Omega \to R$   is a bounded function.
  Then for any
  initial function $\varphi\in C,$ initial value problem (\ref{sdd8-1}), (\ref{sdd8-ic})
  has a global mild solution which satisfies $u\in C([-r,+\infty); L^2(\Omega))$.
  }%


The existence of a mild solution is a consequence of the
continuity of $F(\varphi)\equiv\int_\Omega b(\varphi (-\eta
(\varphi), y)) f(\cdot -y) dy
: C\to L^2(\Omega)$ (see (\ref{sdd8-1})) which gives the
possibility to use the standard method based on Schauder
fixed point theorem (see e.g. \cite[theorem 2.1,
p.46]{Wu_book}). The solution is also global (is defined
for all $t\ge -r$) see e.g. \cite[theorem 2.3,
p.49]{Wu_book}.

\medskip

{\bf Remark~1.} It is important to notice that even in the case of
ordinary differential equations (even scalar) the mapping of the form
$\widetilde{F} (\varphi)=f(\varphi(-r(\varphi)))\, :\,  C([-r_0,0];R) \to
R$ has a very unpleased property. The authors in
\cite[p.3]{Louihi-Hbid-Arino-JDE-2002} write "Notice that the functional
$\widetilde{F}$ is defined on $C,$ but it is clear that it is neither
differentiable nor locally Lipschitz continuous, whatever the smoothness
of $f$ and $r.$" As a consequence, the Cauchy problem associated with
equations with such a nonlinearity "...is {{\texttt{not}} well posed in
the space of continuous functions, due to the non-uniqueness of solutions
whatever the regularity of the functions $f$ and $r$"
\cite[p.2]{Louihi-Hbid-Arino-JDE-2002}. See also a detailed discussion in
\cite{Hartung-Krisztin-Walther-Wu-2006}.

\bigskip

{\bf 4. Uniqueness and well-posedness} \medskip


As in the previous section, we assume that $\eta : C \to
R_{+}$ is continuous and   $f: \Omega -\Omega \to R$   is a
bounded function ($|f(\cdot )|\le M_f$). Additionally, we
assume the following assumption on the delay function
$\eta$ is satisfied

\begin{itemize}
 \item $\exists \eta_{ign}>0$ such that $\eta$ "ignores" values of
 $\varphi(\theta)$ for $\theta\in (-\eta_{ign},0]$ i.e. 
 $$\hskip-10mm \exists\, \eta_{ign}>0 : 
 \forall\varphi^1, \varphi^2\in C :
 \forall\theta\in
 [-r,-\eta_{ign}],\,\Rightarrow \varphi^1(\theta)= \varphi^2(\theta)\quad
   \Longrightarrow \quad
 \eta (\varphi^1)=\eta (\varphi^2). \eqno(H) $$
\end{itemize}
\medskip

{\bf Remark~2.} {\it It is easy to present many examples of (delay)
functions $\eta$, which satisfy assumption (H) . Some of them are
$$\eta(\varphi)= p_1(\varphi(-r)) \hbox{ with } p_1: L^2(\Omega)\to R_{+};$$
$$\eta(\varphi)= \sum^N_{k=1} p_k(\varphi(-r_k)) \hbox{ with } p_k: L^2(\Omega)\to R_{+};
\quad \min r_k >0;$$
$$\eta(\varphi)= \int^{-\eta_{ign}}_{-r} p_1(\varphi(\theta))\, d\theta,\quad
\hbox{ and } \quad\eta(\varphi)= p_1\left( \int^{-\eta_{ign}}_{-r}
\varphi(\theta)\,
 d\theta\right),\quad \eta_{ign}>0,  \hbox{ etc. } $$
 }%

\medskip

{\bf Remark~3.} {\it It is interesting to notice that in
the case $\eta_{ign}>r,$ we have that the delay function
$\eta$ ignores {\tt all} values of $\varphi(\theta),
\forall \theta\in [-r,0],$ so $\eta (\varphi)\equiv const,
\forall\varphi\in C$ i.e. equation (\ref{sdd8-1}) becomes
equation with constant delay. On the other hand,
assumption, similar to (H) with $\eta_{ign}=0,$ is trivial
since $\varphi^1(\theta)=\varphi^2(\theta)$ for all
$\theta\in [-r,0]$ means $\varphi^1=\varphi^2$ in $C,$ so
$\eta (\varphi^1)=\eta (\varphi^2).$ }

 \medskip

To show that 
assumption (H)  implies the {\it uniqueness} of mild solutions (given by
Proposition~1), we will use the standard method of steps with a step less
than 
$\eta_{ign}>0.$ First, let us introduce, for any $\varphi\in C$
 the extension function
$$
\overline \varphi(s)\equiv \left[\begin{array}{ll}
  \varphi(s) & s \in [-r, 0]; \\
  \varphi(0) & s\in (0, \eta_{ign}) \\
\end{array}.
\right. $$

Consider any mild solution $u(t)$ of IVP (\ref{sdd8-1}),
(\ref{sdd8-ic}) and the nonlinearity $\int_\Omega
b(u(t-\eta (u_t), y)) f(\cdot -y) dy.$
 For all $t\in [0,\eta_{ign})$
assumption (H)  gives $\eta(u_t)=\eta(\overline
\varphi_t).$
Let us denote by $r^\varphi(t)\equiv \eta(\overline \varphi_t), t\in
[0,\eta_{ign}).$

Hence any mild solution $u(t)$ of IVP (\ref{sdd8-1}), (\ref{sdd8-ic}) for
all values of $t\in [0, \eta_{ign})$ is also a solution of
\begin{equation}\label{sdd8-4-01}
\left\{\begin{array}{l}
 \dot u(t)+ Au(t) + d\cdot u(t)=
\int_\Omega b(u(t-r^\varphi(t), y)) f(\cdot -y) dy, \quad
t\in [0,\eta_{ign}),\\
  u(\theta)=\varphi(\theta), \, \theta \in  [-r,0],  \\
\end{array}
\right.
\end{equation}
 where $ r^\varphi(t)$ is {\it time}-dependent (but not {\it
state}-dependent delay) known function for all $t\in [0, \eta_{ign}).$ To
show that the last Cauchy problem (with time-dependent delay) has the
unique solution, it is sufficient to consider any two solutions of
(\ref{sdd8-4-01}) and their difference $w(t)\equiv u^1(t)-u^2(t),$ which
satisfies (c.f. (\ref{sdd8-3-1}))
 $$
w(t)= \int^{t}_0 e^{- A (t-s)} \times
$$\begin{equation}\label{sdd8-4-02}\left\{ \int_\Omega
\left[ b(u^1(s-r^\varphi(s), y)) -b(u^2(s-r^\varphi(s), y))
\right] f(\cdot -y) dy
 - d \cdot
w(s)\right\}\, ds, \, t\in [0, \eta_{ign}).
 \end{equation}

An easy calculation, the local Lipschitz property of $b$
and $||e^{- A (t-s)}||\le 1$ give $$ ||w(t)|| \le
\int^{t}_0 (M_f |\Omega| L_b+d)\max_{s\in [0, t]} || w(s)
||\, ds \le t \cdot (M_f |\Omega| L_b+d)  \max_{s\in [0,
t]} || w(s) ||. $$ Here we denote $|\Omega|\equiv
\int_\Omega 1\cdot dx.$

\medskip

{\bf Remark~4.} {\it Here we used that $\max_{s\in [-r, t]}
|| w(s) ||=\max_{s\in [0, t]} || w(s) ||$ since $w(s)\equiv
0$ for $\theta\in [-r,0]$ ($w$ is the difference of two
solutions, both satisfying
(\ref{sdd8-ic})).} %

Choose small enough $\alpha>0 $ to satisfy  $\alpha  (M_f
|\Omega| L_b+d) <1$. The last estimate gives $$ \max_{s\in
[0, \alpha]} || w(s) || \le \alpha \cdot  (M_f |\Omega|
L_b+d) \max_{s\in [0, \alpha]}
 || w(s)
|| < \max_{s\in [0, \alpha]} || w(s) || $$ which implies $\max_{s\in [0,
\alpha]} || w(s) ||=0.$
 This means  that any two mild solutions of  (\ref{sdd8-4-01})  coincide
 for $t\in [0, \alpha]$ with $\alpha < (M_f |\Omega| L_b+d)^{-1}.$ Repeat
  this
 considerations (if necessary) by steps of length $\alpha$ to cover $[0, \eta_{ign}
 ).$ This gives the uniqueness of solutions of (\ref{sdd8-4-01}) and hence the
uniqueness of solutions of (\ref{sdd8-1}), (\ref{sdd8-ic}) for $ t\in [0,
\eta_{ign}).$ The uniqueness on any interval $[0,T]$ holds by the similar
arguments on each time interval $[k\cdot \eta_{ign},
(k+1)\cdot\eta_{ign}), k\in N$ (redefining the function $r^\varphi(t)$).

\smallskip

We may define the {\tt evolution operator} $S_t: C\to C$ by the formula
$S_t \varphi \equiv u_t,$ where $u(t)$ is the unique mild solution of
(\ref{sdd8-1}), (\ref{sdd8-ic}) with initial function $\varphi.$

\medskip

{\bf Remark~5.} {\it The system becomes much simpler if we
additionally assume that the delay function $\eta$
satisfies $$\exists\,  \eta_{min}>0 \, \hbox{ such that }
\, \eta : C \to [\eta_{min}, r] \, \hbox{ is continuous. }
\, \quad \eqno(H1)$$ In that case we may use the classical
method of steps with a step less than $\min \{ \eta_{ign},
\eta_{min} \}.$ To satisfy assumption (H1) for the
functions given in Remark~2 it is sufficient to assume that
$\inf p_i(\cdot) >0. $
}%
\medskip

{\bf Remark~6. } {\it For applied problems described by
ordinary differential equations, condition $\eta (\cdot)
\in [\eta_{min}, r]$ (see (H1)) is used and motivated in
\cite[p.15]{Omari-Gourley} and \cite{Aiello-Freedmand-Wu}.
The authors write "This assumption is known to be realistic
in the case of Antarctic whale and seal populations." }
\medskip

{\bf Remark~7.} {\it It is interesting to notice that
without (H1), one cannot say that the nonlinearity $F(u_t)$
depends on $u|_{[-r,-\eta_{ign}]}$ only. The SD-delay may
vanish (i.e. (H1) does not hold).}

\bigskip

The main result of this section is the following

\medskip

{\bf Theorem~1.} {\it Assume the function $b : R\to R$ is a
locally Lipschitz map, satisfying $|b(w)|\le C_1 |w|+C_b$ with $C_i\ge 0$, 
the delay function $\eta : C\to R_{+}$ is continuous and
satisfies the assumption (H), $f: \Omega -\Omega \to R$ is
a bounded function ($|f(z )|\le M_f, \forall z\in \Omega
-\Omega$). Then the pair $(S_t, C)$ constitutes a dynamical
system i.e. the following properties are satisfied:
\begin{enumerate}
 \item $S_0=Id$ ( identity operator in $C$ );
 \item $\forall\,\, t,\tau \ge 0\quad  \Longrightarrow \quad  S_t\, S_\tau = S_{t+\tau}$;
 \item $t\mapsto S_t$ is a strongly continuous in $C$ mapping;
 \item for any $t\ge 0$ the evolution operator $S_t$ is continuous in $C$ i.e. for any
 $\{\varphi^n\}^\infty_{n=1}\subset C$ such that $||\varphi^n -\varphi||_C\to 0$ as
 $n\to \infty$, one has $||S_t\varphi^n -S_t\varphi||_C\to 0$ as
 $n\to \infty.$
\end{enumerate}
}
\medskip

{\bf Remark~8.} {\it Theorem~1 particularly means that the initial
value problem (\ref{sdd8-1}), (\ref{sdd8-ic}) is {\tt well posed}
in the space $C$ in the sense of J.~Hadamard
\cite{Hadamard-1902,Hadamard-1932}.}

\medskip

{\bf Remark~9.} {\it It is important to emphasize, that we
do not assume the SD-delay function $\eta$ to be Lipschitz
(c.f. \cite{Rezounenko_JMAA-2007}). We propose an
alternative approach, based on the assumption (H) which is
of different nature to the
Lipschitz property of $\eta$. }%

\medskip

 {\it Proof.}
Properties 1), 2) are consequences of the uniqueness of mild solutions due
to (H)  (see considerations above). Property 3) is given by Proposition~1
since the solution is a continuous function $u\in C([-r,T]; L^2(\Omega)).$

Let us prove property 4. Let us fix $t_0\in [0, \eta_{ign}).$ Denote by
$u^k(t)$ the solution of (\ref{sdd8-1}),(\ref{sdd8-ic}) with the initial
function $\varphi^k$ and by $u(t)$ the solution of
(\ref{sdd8-1}),(\ref{sdd8-ic}) with the initial function $\varphi.$

We use the variation of constants formula for parabolic equation
(\ref{sdd8-4-01}) (with $\widetilde A \equiv A+d\cdot$)

\begin{equation}\label{sdd8-4-1}
u(t)=e^{-\widetilde A t}u(0) + \int^{t}_0 e^{-\widetilde A
(t-\tau)} \int_\Omega b(u(\tau-\eta (u_\tau), y)) f(\cdot
-y) dy
\, d\tau,
\end{equation}
\begin{equation}\label{sdd8-4-2}
u^k(t)=e^{-\widetilde A t}u^k((0) + \int^{t}_0 e^{-\widetilde A (t-\tau)}
\int_\Omega b(u^k(\tau-\eta (u^k_\tau), y)) f(\cdot -y) dy
\, d\tau.
\end{equation}
Using $||e^{-\widetilde A t}||\le 1$ and  $||e^{-\widetilde A
(t-\tau)}||\le 1$, we get

$$||u^k(t)-u(t) ||\le ||u^k(0)-u(0)|| $$ $$+ \int^{t}_0
||\int_\Omega \left[ b(u^k(\tau-\eta (u^k_\tau), y)) -
b(u(\tau-\eta (u_\tau), y)) \right] f(\cdot -y) dy|| \,
d\tau $$
\begin{equation}\label{sdd8-4-3}
= ||\varphi^k(0)-\varphi(0)|| + {\cal J}^k_1(t) + {\cal J}^k_2(t),
\end{equation}
where we denote
\begin{equation}\label{sdd8-4-4}
{\cal J}^k_1(s)\equiv \int^{s}_0 || \int_\Omega \left[
b(u^k(\tau-\eta (u^k_\tau), y)) - b(u(\tau-\eta (u^k_\tau),
y)) \right] f(\cdot -y) dy
|| \, d\tau,
\end{equation}
\begin{equation}\label{sdd8-4-5}
{\cal J}^k_2(s)\equiv \int^{s}_0 || \int_\Omega \left[
b(u(\tau-\eta (u^k_\tau), y)) - b(u(\tau-\eta (u_\tau), y))
\right] f(\cdot -y) dy
|| \, d\tau.
\end{equation}
Using the Lipschitz property of $b$, one easily gets 
\begin{equation}\label{sdd8-4-6}
{\cal J}^k_1(t)\le M_f |\Omega| L_b \int^{t}_0 ||u^k(\tau-
\eta(u^k_\tau))- u(\tau- \eta(u^k_\tau))|| \, d\tau \le M_f
|\Omega| L_b t \max_{s\in [-r,t]} ||u^k(s)- u(s)||.
\end{equation}
Estimates (\ref{sdd8-4-6}), (\ref{sdd8-4-3}) and property
${\cal J}^k_2(s) \le {\cal J}^k_2(t)$ for $s\le t \le
t_0<\eta_{ign}$ give $$\max_{t\in [0,t_0]} ||u^k(t)- u(t)||
\le ||\varphi^k(0)-\varphi(0)|| + M_f |\Omega| L_b t_0
\max_{s\in [-r,t_0]} ||u^k(s)- u(s)|| + {\cal J}^k_2(t_0).
$$ Hence $$\max_{s\in [-r,t_0]} ||u^k(s)- u(s)|| \le
||\varphi^k-\varphi||_C + M_f |\Omega| L_b t_0 \max_{s\in
[-r,t_0]} ||u^k(s)- u(s)|| + {\cal J}^k_2(t_0). $$

We choose $t_0< [M_f |\Omega| L_b]^{-1}$ (to satisfy $M_f
|\Omega| L_b t_0 <1$) and get
\begin{equation}\label{sdd8-4-6a}
(1-M_f |\Omega| L_b t_0) \max_{s\in [-r,t_0]} ||u^k(s)-
u(s)|| \le ||\varphi^k-\varphi||_C  +  {\cal J}^k_2(t_0).
\end{equation}

Our goal is to show that ${\cal J}^k_2(t_0) \to 0$ as $k\to \infty.$ The
Lipschitz  property  of $b$ implies
\begin{equation}\label{sdd8-4-7}
{\cal J}^k_2(t_0)\le M_f |\Omega|  L_b \int^{t_0}_0 ||
u(\tau-\eta (u^k_\tau)) - u(\tau-\eta (u_\tau))
 || \, d\tau.
\end{equation}
 We use the extension functions
$$
\overline \varphi(s)\equiv \left[\begin{array}{ll}
  \varphi(s) & s \in [-r, 0]; \\
  \varphi(0) & s\in (0, \eta_{ign}) \\
\end{array}
\right. \quad \hbox{ and }\quad  \overline \varphi^k(s)\equiv
\left[\begin{array}{ll}
  \varphi^k(s) & s \in [-r, 0]; \\
  \varphi^k(0) & s\in (0, \eta_{ign}) \\
\end{array}.
\right.
$$
It is easy to see that the convergence $||\varphi^k -\varphi||_C\to 0$
implies $||\overline \varphi^k_\tau -\overline \varphi_\tau||_C\to 0$ for
any $\tau\in [0, \eta_{ign}).$

On the other hand, for any $\tau\in [0, \eta_{ign}]$ and any $\theta\in
[-r, -\eta_{ign}]$ we have $u^k_\tau(\theta)=\varphi ^k(\tau+\theta)$,
hence assumption (H)  gives $\eta (u^k_\tau)=\eta (\overline
\varphi^k_\tau)$ for any $\tau\in [0, \eta_{ign}).$ The same arguments
give $\eta (u_\tau)=\eta (\overline \varphi_\tau)$ for any $\tau\in [0,
\eta_{ign}).$

The considerations above show that the convergence $||\varphi^k
-\varphi||_C\to 0$ implies $ |\eta (u^k_\tau)-\eta (u_\tau)|= |\eta
(\overline \varphi^k_\tau)-\eta (\overline \varphi_\tau)| \to 0$ for all
$\tau\in [0, \eta_{ign}).$ Here we used the continuity of $\eta : C\to
R_{+}.$

The last property gives that for all $\tau\in [0, \eta_{ign})$ one has
$\tau - \eta (u^k_\tau) \to \tau - \eta (u_\tau),$ when $k\to \infty.$
Hence the
continuity of the mild solution $u$ 
(the strong continuity in $L^2(\Omega)$) implies
(see the integral in (\ref{sdd8-4-7}))
\begin{equation}\label{sdd8-4-8}
\forall \, \tau\in [0, t_0] \Longrightarrow ||u(\tau- \eta(u^k_\tau))-
u(\tau- \eta(u_\tau))|| \to 0, \quad \hbox{ when } \quad k\to \infty.
\end{equation}
On the other hand, it is evidently that
\begin{equation}\label{sdd8-4-9}
\forall \, \tau\in [0, t_0] \Longrightarrow ||u(\tau-
\eta(u^k_\tau))- u(\tau- \eta(u_\tau))|| \le 2\max_{s\in [-r,t_0]}
||u(s)||<+\infty.
\end{equation}
Properties (\ref{sdd8-4-8}) and (\ref{sdd8-4-9}) allow us to use
Lebesgue-Fatou lemma (\cite[p.32]{yosida}) for the scalar function
$||u(\tau- \eta(u^k_\tau))- u(\tau- \eta(u_\tau))||$ to conclude that
\begin{equation}\label{sdd8-4-10}\int^{t_0}_0 ||u(\tau-
\eta(u^k_\tau))- u(\tau- \eta(u_\tau))|| \, d\tau \to 0 \quad \hbox{ when
} \quad k\to \infty.
\end{equation}
Estimates (\ref{sdd8-4-10}) and (\ref{sdd8-4-7}) prove that ${\cal
J}^k_2(t_0) \to 0$ as $k\to \infty.$




Since
$$
\max_{t\in[0,t_0]}{\cal J}^k_2(t)\le {\cal J}^k_2(t_0) \to 0\quad \hbox{
as } k\to \infty,
$$
we finally conclude (see (\ref{sdd8-4-6a}) and the last properties) that
for all $t\in [0,t_0]:$
$$
||u^k_t - u_t||_C\equiv \max_{\theta\in [-r,0]} ||u^k(t+\theta) -
u(t+\theta)||_{L^2(\Omega)} $$
\begin{equation}\label{sdd8-4-11}
\le  [1- M_f |\Omega| L_b t_0]^{-1}\cdot
\left(||\varphi^k-\varphi||_C + {\cal J}^k_2(t_0)\right)
\to 0 \quad \hbox{ as } k\to \infty.
\end{equation}

Estimate (\ref{sdd8-4-11}) gives that $||u^k_t - u_t||_C\to 0$ when
$||\varphi^k -\varphi||_C\to 0$ for all
$t\in [0, t_0] \subset [0, \min \{ [M_f |\Omega| L_b]^{-1}, \eta_{ign}\}).$ 
This is the strong continuity in the space $C$ of the evolution operator
$S_t$ for all (small)
$t\in [0, t_0] \subset [0, \min \{ [M_f |\Omega| L_b]^{-1}, \eta_{ign}\}).$ 

For any $t\ge 0$ we present $S_t\varphi$ as the composition of mappings
\begin{equation}\label{sdd8-4-12}
S_t\varphi = \underbrace{ %
S_p \circ S_p \circ S_p \circ \ldots \circ S_p }_{ q\,\, times } \circ
S_{t-\left[ t\cdot ( 2p)^{-1}\right]}\varphi,\quad \hbox{ where } \quad
p\equiv {1\over 2} \min \{ L^{-1}_b, \eta_{ign}\}, \quad q\equiv \left[
t\over 2p\right].
\end{equation}
Here $\left[ \cdot\right]$ denotes the integer part of a real number.
 The continuity of $S_t$ follows from the proved continuity
of $S_p$ and $S_{t-\left[ t\cdot ( 2p)^{-1}\right]}$ (since
both $p$ and $ t-\left[ t\cdot ( 2p)^{-1}\right]$ belong to
$[0, \min \{ [M_f |\Omega| L_b]^{-1}, \eta_{ign}\})$). The
property 4 is proved. It complies the proof of Theorem~1.
\finproof

\bigskip

\centerline{\bf 5. Asymptotic behavior}

\medskip

This section is devoted to the study of the long-time asymptotic
behavior of the dynamical system $( S_t,C )$, constructed in
Theorem~1.

\medskip

{\bf Theorem~2.} {\it Assume the function $b : R\to R$ is a
locally Lipschitz bounded map ($|b(w)|\le C_b$ with $C_b\ge
0$) and the delay function $\eta : C\to R_{+}$ is
continuous and satisfies the assumption (H),  $f: \Omega
-\Omega \to R$ is a bounded function ($|f(\cdot )|\le
M_f$).
 Then the dynamical
system $( S_t,C )$ has a compact global attractor which is a
compact set in all spaces $C_\delta\equiv C([-r,0]; D(A^\delta)),
\forall\delta\in [0,{1\over 2}).$}

\medskip

{\it Proof.} Our proof is split on four steps.

{\it Step~1.} Let us first prove that for  any  $T>0$ and any
bounded set $B\subset C$
 there exists a constant $C_{T}(B)$ such that
 for any mild solution of (\ref{sdd8-1}), (\ref{sdd8-ic}) with initial values in $B$,
 one has
\begin{equation}\label{sdd8-5-1}
\forall\,  T>0, \forall B, \, \exists\, C_{T}(B)>0 : \forall t\in
[0,T] \Longrightarrow ||u(t)|| \le C_{T}(B).
\end{equation}

Equation (\ref{sdd8-3-1}) implies $||u(t)||\le
||\varphi(0)|| + \int^t_0 (||F(u_s)|| + d ||u(s)||)\, ds. $
Using $||F(u_s)|| \le M_f |\Omega|^{3/2} C_b,$ we have
$||u(t)||\le ||\varphi(0)|| + d \int^t_0 ||u(s)||\, ds + t
M_f |\Omega|^{3/2} C_b.$ Denote by $\Psi (t)\equiv \int^t_0
||u(s)||\, ds$ and use Gronwall lemma to get $\Psi(t)\le
||\varphi(0)|| d^{-1} e^{dt} + e^{dt} M_f |\Omega|^{3/2}
C_b d^{-2} [1-e^{-dt}(dt-1)].$ The last estimate gives
$$||u(t)||\le ||\varphi(0)|| + d\cdot \left[ ||\varphi(0)||
d^{-1} e^{dt} + e^{dt} M_f |\Omega|^{3/2} C_b d^{-2}
[1-e^{-dt}(dt-1)]\right] +  t M_f |\Omega|^{3/2} C_b.$$ It
implies (\ref{sdd8-5-1}).

\medskip

 {\it Step~2.} Next we show that a solution becomes more smooth
 for positive $t.$

 {\bf Lemma.} {\it For any
$0<\alpha<1,\, \epsilon
>0,$
 $T>0$ and any bounded set
 $B\subset C$
 there exists a constant $C_{\alpha,\epsilon,T}(B)$ such that
 for any mild solution of (\ref{sdd8-1}), (\ref{sdd8-ic}) with initial values in $B$,
 one has
\begin{equation}\label{sdd-ls}
 \Vert A^\alpha u(t)\Vert \le C_{\alpha,\epsilon,T}(B)
 \quad \hbox{ for }\,\, t\in (\epsilon,T].
\end{equation}
}%

The proof of the lemma is standard (see e.g. \cite{Chueshov_book}
and also {Rezounenko-MAG-1997, Rezounenko-Wu-2006}).

\medskip

 {\it Step~3.} {\it Dissipativness.} Lemma implies that
 $u(t)\in D(A^{{1\over 2}+\delta}), \forall\delta\in [0,{1/ 2}), t\ge
 \varepsilon(\delta).$ Hence $Au(t)\in D(A^{-{1\over 2}+\delta}), \forall\delta\in [0,{1/ 2}), t\ge
 \varepsilon(\delta).$ Equation (\ref{sdd8-1}) gives $\dot u(t)\in D(A^{-{1\over 2}+\delta}), \forall\delta\in [0,{1/ 2}), t\ge
 \varepsilon(\delta).$ So for all $t\ge \varepsilon(\delta)$ one
 has
 $$ {1\over 2} {d\over dt} || A^\delta u(t)||^2 + || A^{{1\over 2}+\delta}
 u(t)||^2 + d\cdot || A^\delta u(t)||^2  \le ||F(u_t)||\cdot || A^{2\delta}
 u(t)||
 $$
$$\le M_f |\Omega|^{3/2} C_b\cdot \lambda^{2\delta-1}_1||
A^{{1\over 2}+\delta}
 u(t)|| \le {1\over 2} M_f^2 |\Omega|^{3} C^2_b \lambda^{4\delta-2}_1 + {1\over
 2}|| A^{{1\over 2}+\delta}
 u(t)||^2.
$$
Hence $  {d\over dt} || A^\delta u(t)||^2 + || A^{{1\over
2}+\delta} u(t)||^2 + 2d\cdot || A^\delta u(t)||^2  \le
 M_f^2 |\Omega|^{3} C^2_b  \lambda^{4\delta-2}_1. $ Using $|| A^{{1\over
2}+\delta} v||^2 \ge \lambda_1 || A^{\delta} v||^2 $ we have $
{d\over dt} || A^\delta u(t)||^2 + (\lambda_1+2d)|| A^{\delta}
u(t)||^2
  \le  M_f^2 |\Omega|^{3} C^2_b  \lambda^{4\delta-2}_1.
$ Gronwall's  lemma gives $|| A^\delta u(t)||^2 \le ||
A^\delta u(\varepsilon)||^2\cdot \exp\{-(\lambda_1+2d)t\} +
   M_f^2 |\Omega|^{3} C^2_b  \lambda^{4\delta-2}_1(\lambda_1+2d)^{-1}.
$ By lemma (step~2), the value $|| A^\delta u(\varepsilon)||$ is
finite,
which implies that 
$$\forall\delta\in [0,{1/ 2}) \quad \exists\, C(\delta)>0 : \forall B -
\hbox{ bounded in }C, \exists t(B,\delta) :$$
\begin{equation}\label{sdd8-5-3}
\forall t\ge t(B,\delta) \Longrightarrow ||A^\delta u(t)|| \le
C(\delta).
\end{equation}

\medskip

{\it Step~4.} Our next step is to show that the set $\{ S_t\varphi
\, | \, \varphi\in B,\, t>r+\varepsilon \}$ is an equicontinuous
family in $C_\delta\equiv C([-r,0]; D(A^\delta)), \forall\delta\in
[0,{1\over 2}).$

\smallskip

{\bf Remark~10.} {\it In our case we cannot use \cite[thm.
1.8, p.42 ]{Wu_book} since our nonlinearity $F$ is not
Lipschitz.}

\smallskip

We denote by $\{ e_k\}^\infty_{k=1}$  the orthonormal basis (of
$L^2(\Omega)$) of eigenvectors of the operator $A$, so $Ae_k
=\lambda_k e_k,$ $0<\lambda_1<\lambda_2<\ldots <\lambda_k\to
+\infty.$
 Consider for $v\in L^2(\Omega)$
\begin{equation}\label{sdd8-5-4}
\Vert A^\delta \left( e^{-At_1}v - e^{-At_2}v \right)\Vert^2 =
 \sum^\infty_{k=1} \left( e^{-\lambda_kt_1} - e^{-\lambda_kt_2}
 \right)^2 \lambda^{2\delta}_k \cdot v^2_k,
 \quad \hbox{ where }\,\, v_k\equiv \langle v, e_k\rangle.
\end{equation}

Assuming $0<t_1<t_2,$ one can easily check that $
\left(
e^{-\mu t_1} - e^{-\mu t_2}
 \right)^2 \mu^{2\delta} = \left|  e^{-\mu t_1} - e^{-\mu t_2}
  \right| \cdot \left| \left( e^{-\mu t_1} - e^{-\mu t_2}
 \right) \mu^{2\delta} \right|  \le
\left| \left( e^{-\mu t_1} - e^{-\mu t_2}
 \right) \mu^{2\delta} \right| \le |t_1-t_2|\cdot\max_{\tau\in [t_1,t_2]}
 \mu^{1+2\delta}e^{-\mu \tau}
 =|t_1-t_2|\cdot \mu^{1+2\delta}e^{-\mu t_1}
 ,\, \mu>0.$ Calculations give $\max_{\mu>0} \mu^{1+2\delta}e^{-\mu
 t_1}=e^{-(1+2\delta)} \left({1+2\delta\over
 t_1}\right)^{1+2\delta}.$ Hence, for any $k\in N$, one has
 $$\left( e^{-\lambda_kt_1} - e^{-\lambda_kt_2}
 \right)^2 \lambda^{2\delta}_k \le |t_1-t_2| \cdot e^{-(1+2\delta)}
 (1+2\delta)^{(1+2\delta)} \cdot \left(\min\{ t_1,t_2\}\right)^{-(1+2\delta)}
 $$
 The last estimate and (\ref{sdd8-5-4}) give
\begin{equation}\label{sdd8-5-5}
\Vert A^\delta \left( e^{-At_1}v - e^{-At_2}v \right)\Vert \le
D_\delta \left(\min\{ t_1,t_2\}\right)^{-\left( \delta+{1\over
2}\right) }\cdot \sqrt{|t_1-t_2|}\cdot ||v||,
\end{equation}
where $D_\delta \equiv e^{-\left( \delta+{1\over 2}\right) } \left(
\delta+{1\over 2}\right)^{\left( \delta+{1\over 2}\right) }.$ In the
same way we get
\begin{equation}\label{sdd8-5-5old}
\Vert A^\delta \left( e^{-At_1}v - e^{-At_2}v \right)\Vert \le
\widehat D_\delta \left(\min\{ t_1,t_2\}\right)^{-(1+\delta)}\cdot
|t_1-t_2|\cdot ||v||,
\end{equation}
where $\widehat D_\delta \equiv e^{-(1+\delta)}
(1+\delta)^{1+\delta}$ and also (see e.g.
\cite[(2.8.6)]{Chueshov_book}), we get $|| A^\delta e^{-As}||\le
(e\cdot s)^{-\delta} \delta^\delta.$ Calculations and the last
estimate give for $0<t_1<t_2$
\begin{equation}\label{sdd8-5-5a}
\int^{t_2}_{t_1} \Vert A^\delta  e^{-A(t_2-\tau)}\Vert \,
d\tau \le \left(\frac{\delta}{e}\right)^\delta \cdot
\frac{|t_1-t_2|^{1-\delta}}{1-\delta}.
\end{equation}

Let us consider for any mild solution $u(t)$, the function
$G(t)\equiv F(u_t) + d u(t)$ and the difference (see
(\ref{sdd8-3-1})) $$ \Vert A^\delta \left( u(t_1) - u(t_2)
\right)\Vert \le \Vert A^\delta \left( e^{-At_1}\varphi(0) -
e^{-At_2}\varphi(0) \right)\Vert $$ $$+ \int\limits_0^{t_1} \Vert
A^\delta \left( e^{-A(t_1-\tau)} G(\tau) - e^{-A(t_2-\tau)} G(\tau)
\right)\Vert\, d\tau + \int\limits_{t_1}^{t_2} \Vert A^\delta
e^{-A(t_2-\tau)} G(\tau)\Vert\, d\tau. $$ Here, as before
$0<t_1<t_2.$ The last estimate, (\ref{sdd8-5-5}) and
(\ref{sdd8-5-5a}) give
\begin{equation}\label{sdd8-5-6}
\Vert A^\delta \left( u(t_1) - u(t_2) \right)\Vert \le
L(\delta,t_1,t_2,\varphi) \cdot \sqrt{|t_1-t_2|},
\end{equation}
where
$$
 L(\delta,t_1,t_2,\varphi) \equiv D_\delta
\left(\min\{ t_1,t_2\}\right)^{-\left( \delta+{1\over 2}\right)}
||\varphi(0)||
$$
\begin{equation}\label{sdd8-5-7}
+\left[ D_\delta t_1^{{1\over 2}-\delta} \left({1\over
2}-\delta\right)^{-1} + \delta^\delta
\left(e(1-\delta)\right)^{-1}\right] \cdot \max_{\tau\in [0,t_2]}
||F(u_\tau) + d u(\tau)||.
\end{equation}

Here we used
$$\int\limits_0^{t_1} \Vert
A^\delta \left( e^{-A(t_1-\tau)} G(\tau) - e^{-A(t_2-\tau)} G(\tau)
\right)\Vert\, d\tau $$ $$\le D_\delta \sqrt{|t_1-t_2|}
\max_{\tau\in [0,t_2]} ||G(\tau)|| \cdot \int\limits_0^{t_1}
(t_1-\tau)^{-\left( \delta+{1\over 2}\right)}\, d\tau $$ $$ =
D_\delta \sqrt{|t_1-t_2|} \max_{\tau\in [0,t_2]} ||G(\tau)|| \cdot
t_1^{{1\over 2}-\delta} \left({1\over 2}-\delta\right)^{-1}.$$

It is evidently that 
(\ref{sdd8-5-3}) gives $\max_{\tau\in [0,t_2]} ||F(u_\tau) + d
u(\tau)|| <C$ for any mild solution which is already in the ball of
dissipation i.e. (\ref{sdd8-5-3}) holds. These, the form of the
constant $ L(\delta,t_1,t_2,\varphi)$ (see (\ref{sdd8-5-7})) and
(\ref{sdd8-5-6}) imply that for any time interval $[a,b]\subset
(0,+\infty)$ with $a> t(B,\delta)$ (see (\ref{sdd8-5-3})), there
exists a constant $L>0$ such that for any mild solution $u(t)$ and
any $t_1,t_2\in [a,b]$ one has
\begin{equation}\label{sdd8-5-8}
\Vert A^\delta \left( u(t_1) - u(t_2) \right)\Vert \le L \cdot
\sqrt{|t_1-t_2|}.
\end{equation}
 This gives the equicontinuity of the
family $\{ u(t+\theta)\, | \, \theta\in [-r,0], t> t(B,\delta) \}$
in all spaces $C_\delta\equiv C([-r,0]; D(A^\delta)),
\forall\delta\in [0,{1\over 2}).$

\medskip

Finally,  estimate (\ref{sdd8-5-8}) for $\delta\in [0,{1\over 2})$
and estimate (\ref{sdd8-5-3}) for $\delta_1\in (\delta,{1\over 2})$,
particularly mean (by Arzela-Ascoli theorem) that for any
$\varphi\in C$ and $t>t(B,\delta)$ (see (\ref{sdd8-5-3}), one has
$S_t\varphi\in K_\delta, $  where $K_\delta$ is a compact set in all
spaces $C_\delta\equiv C([-r,0]; D(A^\delta)), \forall\delta\in
[0,{1\over 2}).$

\medskip

That means that the dynamical system $(S_t;C)$ is dissipative and
asymptotically compact, hence by the classical theorem on the
existence of an attractor (see, for example,
\cite{Babin-Vishik,Temam_book})
$(S_t;C)$ has a compact global attractor. The proof of Theorem~2
is complete. \rule{5pt}{5pt}
\medskip

{\bf Remark~11. } {\it All the results above are valid for
local in space variable nonlinearity i.e. equation
(\ref{sdd8-1}) with
\begin{equation}\label{sdd8-5-10}
\big( F_\ell (u_t)
\big)(x) \equiv b(u(t-\eta (u_t), x)),\quad x\in \Omega. 
\end{equation}
}

\medskip

As an application we can consider the diffusive Nicholson's
blowflies equation (see e.g. \cite{So-Yang,So-Wu-Yang})
with state-dependent delays. More precisely, we consider
equation (\ref{sdd8-1}) where $-A$ is the Laplace operator
with the Dirichlet boundary conditions, $\Omega\subset
R^{n_0}$ is a bounded domain with a smooth boundary, the
function $f$ can be,
 for example, $ f(s)={1\over \sqrt{4\pi\alpha}}
e^{-s^2/4\alpha}$, as in \cite{So-Wu-Zou} (for the
non-local in space problem) or Dirac delta-function to get
the local problem (see Remark~11),
the nonlinear function $b$ is
given by $b(w)=p\cdot we^{-w}.$ Function $b$ is bounded
, so for any continuous delay function $\eta$, satisfying (H), the
conditions of Theorems~1,2 are valid. As a result, we conclude
that
the initial value problem (\ref{sdd8-1}),(\ref{sdd8-ic}) is well
posed in $C$ and the dynamical system $(S_t,C)$ has a global
attractor (Theorem~2).



\medskip

{\bf Remark~12. } {\it All the considerations above are
obviously valid for O.D.E.s, for example, of the following
form

\begin{equation}\label{sdd8-5-12}
\dot u (t) + A u(t)+ d\cdot u(t) = b(u(t-\eta (u_t))),
\quad u(\cdot)\in R^n, \, d\ge 0.
\end{equation}
One simply needs to substitute $L^2(\Omega )$ by $R^n$ and use
$C\equiv C([-r,0];R^n)$ instead of  $ C([-r,0];L^2(\Omega )).$
Assumptions on the delay function $\eta$ are the same. The function
$b : R^n\to R^n$ is locally Lipschitz continuous and satisfies
$||b(w)||_{R^n} \le C_1 ||w||_{R^n} +C_b$ with $C_1, C_b \ge 0.$
 For the method of steps in the case of O.D.E.s see e.g.
\cite[Proposition~1]{Walther2002}.

} 

\medskip

\noindent {\bf Acknowledgements.}  The author wishes to thank
Professors I.D.~Chueshov and  \\ H.-O.~Walther for useful
discussions of an early version of the manuscript.
\bigskip

%





\vskip2cm




\hfill Kharkiv

\bigskip

\hfill February 4, 2008

\end{document}